\documentclass[graybox, natbib]{svmult}
\bibpunct{(}{)}{;}{a}{}{,} 

\usepackage{mathptmx}       
\usepackage{helvet}         
\usepackage{courier}        
\usepackage{type1cm}        

\usepackage{makeidx}         
\usepackage{graphicx}        
\usepackage{multicol}        
\usepackage[bottom]{footmisc}
\usepackage[normalem]{ulem}	
\usepackage{hyperref}  
\usepackage{soul}   
\usepackage{psfrag}

\makeindex             

\usepackage{amssymb,epsfig}
\usepackage{amsmath}
\usepackage{color}
\usepackage{enumerate}
\usepackage{psfrag}

\def\div{{\rm div \;}}

\def\N{\mathbb{N}}
\def\R{\mathbb{R}}

\def\E{\mathbb{E}} 
\def\P{\mathbb{P}} 

\begin{document}

\title*{Mathematical foundations of Accelerated Molecular Dynamics methods}
\author{Tony Leli\`evre}

\institute{Tony Leli\`evre \at CERMICS, \'Ecole des Ponts, Universit\'e
    Paris-Est, INRIA, 77455 Champs-sur-Marne, France. \email{tony.lelievre@enpc.fr}}
%
%
\maketitle
\abstract{The objective of this review article is to present recent results on the mathematical analysis of the Accelerated Dynamics algorithms introduced by A.F.~Voter in collaboration with D. Perez and M. Sorensen. Using the notion of quasi-stationary distribution, one is able to rigorously justify the fact that the exit event from a metastable state for the Langevin or overdamped Langevin dynamics can be modeled by a kinetic Monte Carlo model. Moreover, under some geometric assumptions, one can prove that this kinetic Monte Carlo model can be parameterized using Eyring-Kramers formulas. These are the building blocks required to analyze the Accelerated Dynamics algorithms, to understand their efficiency and their accuracy, and to improve and generalize these techniques beyond their original scope.}

\section{Introduction}\label{sec:intro}

The objective of this article is to review recent works on the mathematical analysis of the {\em accelerated dynamics techniques} introduced by A.F.~Voter and co-workers from the late nineties up to recently~\citep{voter-97,voter-98,sorensen-voter-00,perez-cubuk-waterland-kaxiras-voter-15}. The objectives of these methods is to efficiently simulate thermostated molecular dynamics trajectories over very large timescales. The mathematical analysis which has been developed gives the set of assumptions underlying these algorithms, and a way to assess their accuracy. This helps to understand their limitations, but also to improve and extend these numerical methods beyond their original scope.

The aim of statistical computational physics and molecular dynamics is to infer from a microscopic model of matter its macroscopic properties. At the microscopic level, the atomic configuration is given by a set of positions of the atoms (or group of atoms), and the basic ingredient is a potential energy function:
\begin{equation}\label{eq:V}
V:\R^d \to \R.
\end{equation}
Here $d$ is the number of degrees of freedom (typically three times the number of atoms) and the function $V$ associates to a given atomic configuration its energy. We will not discuss here how to build such a function $V$. Let us simply mention that this is related to the construction of so-called force fields, which requires a lot of chemical intuition to infer the functional form of the force field, and good parametrization using either experimental data or {\em ab initio} computations to evaluate the electronic structure associated to given positions of the nuclei. Combining optimally these informations in order to improve the predictive abilities of forces fields is a very lively research subject, at the interface between physics, chemistry, numerical analysis and data sciences.

For a given potential $V$, at a fixed temperature $T$, the configurations of the atomic system at thermal equilibrium are distributed according to the Boltzmann Gibbs measure (a.k.a. the canonical measure)
\begin{equation}\label{eq:mu}
\mu=Z^{-1} \exp(-\beta V(q)) \, dq
\end{equation}
where $q \in \R^d$ is the positions of the atoms, $\beta^{-1}= k_B T$ ($k_B$ being the Boltzmann constant) and $Z= \int_{\R^d} \exp(-\beta V(q)) \, dq$ is assumed to be finite. Computing averages with respect to $\mu$ gives access to so-called {\em thermodynamic quantities}. Examples include heat capacity, free energy difference,  stress, etc. The focus of this work is rather to discuss numerical methods to compute {\em dynamical quantities}, namely observables which depend on the trajectories of the molecular system: transition times, transition paths, reaction mechanisms, etc. This requires to define a dynamics. The typical dynamics one should have in mind is the
Langevin dynamics:
\begin{equation}\label{eq:Langevin}
\left\{
\begin{aligned}
dq_t & = M^{-1} p_t \, dt \\
dp_t &= -\nabla V(q_t) \, dt - \gamma M^{-1} p_t \, dt + \sqrt{2
  \gamma \beta^{-1}} dW_t
\end{aligned}
\right.
\end{equation}
where $(q_t,p_t) \in \R^{d \times d}$ denotes the positions and momenta of the particles
at time $t \ge 0$,
$M$~is the mass tensor, $\gamma >0$ is a friction parameter and $W_t$ is a
standard $d$-dimensional Brownian motion. Under loose
assumptions on $V$, this dynamics is ergodic with respect to the
phase-space canonical measure $\overline{Z} \exp(-\beta (V(q) +p^T M^{-1} p /2)) \, dq \, dp$, where $
\overline{Z}= \int_{\R^d} \exp(-\beta (V(q) +p^T M^{-1} p /2)) \, dq \, dp$. The marginal in $q$ of this measure is the measure $\mu$ defined above (see~\eqref{eq:mu}). In particular, for any test function $\varphi:\R^d
\times \R^d \to \R$,
\begin{equation}\label{eq:ergodic}
\lim_{t \to \infty} \frac{1}{t} \int_0^t \varphi (q_s) \, ds = \int_{\R^d} \varphi d \mu.
\end{equation}
 In the following, we will also consider the overdamped
Langevin dynamics which is obtained from~\eqref{eq:Langevin} in the
limit $\gamma \to \infty$ or $M \to 0$ (see for example~\citet[Section 2.2.4]{lelievre-rousset-stoltz-book-10}):
\begin{equation}\label{eq:overdamped_Lang}
dX_t=-\nabla V(X_t) \, dt + \sqrt{2 \beta^{-1}} dW_t
\end{equation}
where $X_t \in \R^d$ denotes the positions of the particles.
More precisely, the overdamped Langevin
dynamics is for example derived from the Langevin dynamics in the large friction
limit and using a rescaling in time: assuming $M={\rm Id}$ for simplicity,
in the limit $\gamma
\to \infty$, $(q_{\gamma t})_{t \ge 0}$
converges to $(X_t)_{t \ge 0}$.  Again, under loose
assumptions on $V$, this dynamics is ergodic with respect to the
canonical measure $\mu$. The Langevin and overdamped Langevin dynamics are thus thermostated dynamics: they describe the evolution of the system at a given temperature $T$. Simulating these dynamics over large timescales in order to have access to macroscopic properties of the system (both thermodynamic and dynamical quantities) is the objective of molecular dynamics simulations with applications in many scientific areas: biology, chemistry, material sciences, etc. To obtain accurate results, this requires to simulate stochastic dynamics in large dimension over very large time scales.

The main difficulty when simulating the dynamics~\eqref{eq:Langevin} or~\eqref{eq:overdamped_Lang} in practice is that they are metastable. This means that the trajectory of the positions $(q_t)_{t \ge 0}$ or $(X_t)_{t \ge 0}$ remains trapped for very long periods of time in some regions of the configurational space called {\em metastable states}. This is actually expected from a physical viewpoint: these metastable states typically correspond to some macroscopic conformations of the system. In material sciences for example, one could think of these metastable states as positions of some defects in a crystal, or of some ad-atoms on a surface. In biology, these metastable states can be associated for example with molecular conformations of  a protein. And it is indeed expected that the residence times in these metastable states are much larger than the typical timescale of vibration within the metastable states. In practice, one needs for example to use timesteps of the order of $10^{-15} s$ to discretize~\eqref{eq:Langevin}, while the typical phenomena of interest occur over timescales of the order of $10^{-6} s$ up to seconds or even more ! The numerical counterpart is twofolds: to compute thermodynamic quantities, the convergence of time averages such as~\eqref{eq:ergodic} is very slow; and to evaluate dynamical quantities of interest, namely typically transition paths and transition times between metastable states, one needs to sample rare events, namely the exit from metastable states, and the transition to a new metastable state. The objective of the {\em accelerated dynamics} algorithms is indeed to efficiently sample metastable dynamics in order to have access to some dynamical quantities. Let us make precise that we are here considering numerical methods to sample the whole dynamics from states to states, and not for example the ensemble of reactive paths between two given metastable states, for which other dedicated methods can be used (splitting techniques~\citep{van-erp-moroni-bolhuis-03,allen-warren-ten-wolde-05,cerou-guyader-lelievre-pommier-11}, transition path sampling methods~\citep{dellago-bolhuis-chandler-99}, etc.)

The basic idea of the {\em accelerated dynamics} algorithms is the following: if the stochastic process $(q_t,p_t)_{t \ge 0}$ or $(X_t)_{t \ge 0}$ remains trapped for a sufficiently long time in a metastable state, it forgets the way it entered this state, and this means that the exit event from this metastable state can be modeled by the exit event of a kinetic Monte Carlo model, namely a pure jump Markov process. This will be explained in Section~\ref{sec:kMC_QSD}. We will then show in Section~\ref{sec:EK_HTST} that, in the small temperature regime, it is possible to parameterize the kinetic Monte Carlo model modelling the exit event using Eyring-Kramers laws, which gives the basis for a rigorous foundation of the harmonic transition state theory. From these properties, it is then possible to devise efficient algorithms to simulate metastable dynamics. This is explained in Section~\ref{sec:AD}. The bottom line is that the details of the dynamics within metastable states are not interesting: only the exit events from these states need to be simulated.


\section{Kinetic Monte Carlo models and quasi-stationary distribution}\label{sec:kMC_QSD}

In this section, we consider a set $S \subset \R^d$, which is assumed to be bounded, regular and open. This set is intended to be associated to one state of a kinetic Monte Carlo (kMC) model, and we would like to relate the exit event from $S$ using the Langevin~\eqref{eq:Langevin} or overdamped Langevin~\eqref{eq:overdamped_Lang} dynamics and a kinetic Monte Carlo model. We will first present in Section~\ref{sec:kMC} how the exit event from a state is modeled in a kinetic Monte Carlo model. We will then introduce in Section~\ref{sec:QSD} the notion of quasi-stationary distribution (QSD), which is the basic ingredient to connect the simulation of the exit event from a set $S$ for~\eqref{eq:Langevin} or~\eqref{eq:overdamped_Lang} with a kMC model, as explained in Section~\ref{sec:QSD_kMC}. Finally, we will conclude this section with a discussion on the way to estimate the convergence time to the QSD in Section~\ref{sec:tau_corr}.
 
\subsection{Kinetic Monte Carlo models}\label{sec:kMC}

Kinetic Monte Carlo models (a.k.a.  Markov state models, see~\citet{bowman-pande-noe-14,schuette-sarich-13}) are continuous-time Markov processes with values
in a discrete state space, namely a jump Markov process (see~\citet{voter-05} for a nice introduction to kMC models). They consist of a
collection of states that we can assume to be indexed by integers, and rates $(k_{i,j})_{i \neq j
  \in \N}$ which are associated with transitions between these
states. For a state $i \in \N$, the states $j$ such
that $k_{i,j} \neq 0$ are the neighboring states of $i$ denoted in the
following by
\begin{equation}\label{eq:Ni}
{\mathcal N}_i=\{j \in \N, \, k_{i,j} \neq 0\}.
\end{equation}
 One can thus represent such a jump Markov model as a graph: the states are the vertices,
and an oriented edge between two vertices $i$ and $j$ indicates that
$k_{i,j}\neq 0$.

Starting at time $0$ from a state $Y_0 \in \N$, the model consists in iterating the following two steps
over $n \in \N$: Given $Y_n$,
\begin{itemize}
\item Sample the residence time $T_n$ in $Y_n$ as an exponential random
  variable with parameter $\sum_{j \in {\mathcal N}_{Y_n}} k_{Y_n,j}$:
\begin{equation}\label{eq:T_n}
\forall t \ge 0, \, \P(T_n\ge t | Y_n=i) = \exp\left(- \left[\sum_{j
      \in {\mathcal N}_{i}} k_{i,j}\right] \, t\right).
\end{equation}
\item Sample independently from $T_n$ the next visited state $Y_{n+1}$ starting from $Y_n$
  using the following law
\begin{equation}\label{eq:Y_n}
\forall j  \in {\mathcal N}_{i}, \, \P(Y_{n+1}=j | Y_n=i) = \frac{k_{i,j}}{\sum_{j'
    \in {\mathcal N}_{i}} k_{i,j'}}.
\end{equation}
\end{itemize}
The associated continuous-time process $(Z_t)_{t \ge 0}$
with values in $\N$ defined by:
\begin{equation}\label{eq:Z}
\forall n \ge 0, \, \forall t \in \left[\sum_{m=0}^{n-1} T_m,
\sum_{m=0}^{n} T_m\right), \quad Z_t = Y_n
\end{equation}
(with the convention $\sum_{m=0}^{-1}=0$) is then a (continous-time) jump Markov process.

The exit event from the state $i$ is thus modeled by the couple of random variables $(T,Y)$ where:
\begin{enumerate}
\item $T$ and $Y$ are independent ;
\item $T$ is exponentially distributed with parameter $\sum_{j} k_{i,j}$ ;
\item $Y$ takes the value $j \in \mathcal N_i$ with probability $\frac{k_{i,j}}{\sum_{j'
    \in {\mathcal N}_{i}} k_{i,j'}}$.
\end{enumerate}

The question we would like to address in the following is: when can one use such a model to simulate the exit event from a metastable state for the Markov dynamics~\eqref{eq:Langevin} or~\eqref{eq:overdamped_Lang} ? The cornerstone to make this connection is the notion of quasi-stationary distribution, which will be introduced in Section~\ref{sec:QSD}.

The interest of using a kMC model rather than the Markov dynamics~\eqref{eq:Langevin} or~\eqref{eq:overdamped_Lang} to model the evolution of the system is twofold. From a modelling viewpoint, new insights can be
gained by building such coarse-grained models, that are easier to handle. From a numerical viewpoint,
the hope is to be able to build the jump Markov process from short
simulations of the full-atom dynamics from states to states. Then, once the rates have been defined, it is possible to simulate the system over
much larger timescales than the time horizons attained by standard
molecular dynamics, either by using directly the jump Markov process,
or as a support to accelerate molecular
dynamics (see for example~\citet{voter-97,voter-98,sorensen-voter-00} and Section~\ref{sec:AD} below). It is also
possible to use dedicated algorithms to extract from the graph
associated with the jump Markov process the most important features of
the dynamics (for example quasi-invariant sets and essential
timescales using large deviation theory~\citep{freidlin-wentzell-84}), see for example~\citet{wales-03,cameron-14}.

\subsection{Quasi-stationary distribution}\label{sec:QSD}

In this section, we
focus for simplicity on the overdamped Langevin
dynamics~\eqref{eq:overdamped_Lang}. Generalizations to the Langevin dynamics~\eqref{eq:Langevin} are expected to be true, and are the subject of works under progress (see for example~\citet{nier-18} for a first result in that direction). Let us consider the first exit time from a fixed set $S \subset \R^d$:
$$T_S=\inf \{t \ge 0, \, X_t \not\in S \}.$$
The exit event from $S$ is fully characterized by the couple of random variables:
$$(T_S,X_{T_S}).$$

The basic intuition already mentioned in the introduction is that if the process remains for a sufficiently long time in $S$, it should be possible to model the exit event by the exit event of a kMC model. This naturally leads us to consider the quasi-stationary distribution which is the invariant law for the process conditioned to stay in~$S$.

\begin{definition}
A probability measure $\nu_S$ with support in $S$ is called a quasi-stationary distribution (QSD) for
the Markov process $(X_t)_{t \ge 0}$ if and only if 
$$\forall t > 0, \, \forall A \subset S, \nu_S(A)=\frac{\displaystyle \int_S \P(X_t^x
  \in A, t < T^x_S) \, \nu_S(dx)}{\displaystyle \int_S \P(t < T^x_S) \, \nu_S(dx)}.$$
\end{definition}
In other words, $\nu_S$ is a QSD if when $X_0$ is distributed according
to $\nu_S$, for all
positive~$t$, the law of~$X_t$ conditionally to the fact that $(X_s)_{0
  \le s \le t}$ remains in the state $S$ is still $\nu_S$.

The QSD satisfies three properties which will be crucial in the
following. We refer for example to~\citet{le-bris-lelievre-luskin-perez-12} for a proof of these results and to~\citet{collet-martinez-san-martin-13} for more
general results on QSDs.

\begin{proposition}\label{prop:QSD_1}
Let $(X_t)_{t \ge 0}$ follow the dynamics~\eqref{eq:overdamped_Lang}
with an initial condition $X_0 \in S$. Then, there exists a probability
distribution $\nu_S$ with support in $S$ such that
\begin{equation}\label{eq:QSD_1}
\lim_{t \to \infty} {\mathcal L}(X_t | T_S > t) = \nu_S,
\end{equation}
where for a given $t>0$, ${\mathcal L}(X_t | T_S > t)$ denotes the law of the random variable $X_t$ conditioned to the event $\{T_S>t\}$. The distribution $\nu_S$ is the QSD associated with $S$.
\end{proposition}
A consequence of this proposition is the existence and uniqueness of
the QSD. In some sense, the QSD can thus be seen as the longtime limit of the process
conditioned to stay in the state $S$.  

Let us now give a second property of the QSD.
\begin{proposition}\label{prop:QSD_2}
Let $L=-\nabla V \cdot \nabla + \beta^{-1} \Delta$ be the
infinitesimal generator of $(X_t)_{t \ge 0}$ (satisfying~\eqref{eq:overdamped_Lang}). Let us consider the
first eigenvalue and eigenfunction associated with the adjoint operator
$L^*=\div ( \nabla V + \beta^{-1} \nabla)$ with homogeneous Dirichlet
boundary condition:
\begin{equation}\label{eq:QSD_2}
\left\{
\begin{aligned}
L^* u_1 &= - \lambda_1 u_1 \text{ on $S$},\\
u_1 &=0 \text{ on $\partial S$}.
\end{aligned}
\right.
\end{equation}
The QSD $\nu_S$ associated with $S$ satisfies:
$$d\nu_S = \frac{u_1(x) \, dx}{\displaystyle \int_S u_1(x) \, dx}$$
where $dx$ denotes the Lebesgue measure on $S$.
\end{proposition}
The QSD thus has a density with respect to the Lebesgue measure, which
is nothing but the ground state of the Fokker-Planck operator $L^*$
associated with the dynamics with absorbing boundary conditions. The existence, uniqueness and positivity of this ground state is a standard consequence of the ellipticity of the operator $L^*$ and the boundedness of the set $S$.

Finally, the last property of the QSD concerns the exit event from $S$, when
$X_0$ is distributed according to $\nu_S$.
\begin{proposition}\label{prop:QSD_3}
Let us assume that $X_0$ is distributed according to the QSD $\nu_S$ in $S$. Then the
law of the couple $(T_S,X_{T_S})$ (namely the first exit time and the
first exit point from $S$) is fully characterized by the following properties:
\begin{enumerate}
\item $T_S$ is independent of $X_{T_S}$
\item $T_S$ is exponentially distributed with parameter $\lambda_1$
  (defined in Equation~\eqref{eq:QSD_2} above)
\item The law of $X_{T_S}$ is given by\footnote{Here and in the following, the superscript $\nu_S$ in $\E^{\nu_S}$ of $\P^{\nu_S}$ indicates that the stochastic process starts under the QSD: $X_0 \sim \nu_S$.}: for any bounded measurable
  function $\varphi: \partial S \to \R$,
\begin{equation}\label{eq:QSD_3}
\E^{\nu_S}(\varphi(X_{T_S})) = - \frac{\displaystyle\int_{\partial S}
  \varphi \, \partial_n u_1 \, d \sigma}{\displaystyle\beta \lambda_1 \int_{S}
  u_1(x) \, dx}
\end{equation}
where $\sigma$ denotes the Lebesgue measure on $\partial S$ induced by
the Lebesgue measure in $\R^d$ and the Euclidean scalar product, and
$\partial_n u_1=\nabla u_1 \cdot n$ denotes the outward normal
derivative of $u_1$ on $\partial S$. 
\end{enumerate}
\end{proposition}


Proposition~\ref{prop:QSD_3} explains the interest of the QSD. Indeed, if the
process is distributed according to the QSD in $S$ (namely, from
Proposition~\ref{prop:QSD_1}, if it remained for a sufficiently long
time in $S$), then the exit event from the state $S$ can be modeled using a kinetic Monte Carlo model, since the exit time is exponentially distributed and independent of the exit point. This will be explained  in more detail in the next section.

\begin{remark}[From overdamped Langevin to Langevin]\label{rem:lang_1}
As mentioned above, the existence of the QSD and the convergence of the
conditioned process towards the QSD for the Langevin
process~\eqref{eq:Langevin} requires extra work compared to the
overdamped Langevin process~\eqref{eq:overdamped_Lang}. For results in that direction, we refer to the recent
manuscript~\citet{nier-18}. The main difficulties are twofold: (i) even if $S$ is a bounded set, the associated ensemble in phase-space is not bounded (the velocities are indeed not bounded) and (ii) the Langevin dynamics is not reversible and not elliptic (noise only acts on velocities, not on positions). This implies some difficulties when studying the spectral properties of the infinitesimal generator with absorbing boundary conditions on $\partial S$.
\end{remark}

\subsection{Modeling of the exit event from a metastable state}\label{sec:QSD_kMC}

Using Proposition~\ref{prop:QSD_1}, if the process $(X_t)_{t \ge 0}$ remains sufficiently long in the state~$S$, the random variable $X_t$ is approximately distributed according to the QSD $\nu_S$.
Then, from Proposition~\ref{prop:QSD_3}, the exit event $(T_S,X_{T_S})$ satisfies the basic properties needed to be modeled by the exit event of a kinetic Monte Carlo model: $T_S$ is independent from $X_{T_S}$ and exponentially distributed. 

\begin{figure}[h]
\centering
\psfrag{S}{$S$}
\psfrag{S_1}{$S_1$}
\psfrag{S_2}{$S_2$}
\psfrag{S_3}{$S_3$}
\psfrag{S_4}{$S_4$}
\psfrag{dS_1}{$\partial S_1$}
\psfrag{dS_2}{$\partial S_2$}
\psfrag{dS_3}{$\partial S_3$}
\psfrag{dS_4}{$\partial S_4$}
  \includegraphics[height=5cm]{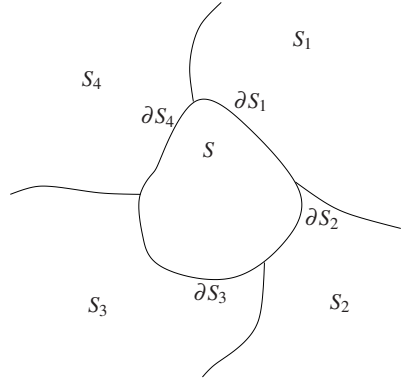}
  \caption{The boundary $\partial S$ of the domain $S$ is divided into
    $4$ subdomains $(\partial S_i)_{1 \le i \le 4}$, which are the
    common boundaries with the neighboring states $(S_i)_{1 \le i \le 4}$.}
  \label{fig:S}
\end{figure}

To be more specific, let us consider the domain $S$, and let us divide its boundary $\partial S$ into subsets $(\partial S_i)_{i=1, \ldots, J}$ associated to transitions to neighboring states $(S_i)_{i=1, \ldots, J}$, see Figure~\ref{fig:S} for a schematic representation in the case $J=4$. The next visited state is thus defined by the random variable $Y_S$ with values in $\{1, \ldots, J\}$ defined by:
$$Y_S=i \text{ if and only if } X_{T_S} \in \partial S_i.$$
Let us now introduce the rates: for $i=1, \ldots, J$,
\begin{equation}\label{eq:ki}
k_i=\lambda_1 \P^{\nu_S}(X_{T_S} \in \partial S_i)
\end{equation}
where from~\eqref{eq:QSD_3}
\begin{equation}\label{eq:pi}
 \P^{\nu_S}(X_{T_S} \in \partial S_i) =
-\frac{\displaystyle \int_{\partial S_i} \partial_n u_1 \, d \sigma}{\displaystyle \beta \lambda_1
  \int_S u_1(x)  \, d x}.
\end{equation}
We recall that $(\lambda_1,u_1)$ has been defined in Proposition~\ref{prop:QSD_2}. Now, from Proposition~\ref{prop:QSD_3}, we obviously have the following properties on the couple of random variables $(T_S,Y_S)$:
\begin{enumerate}
\item $T_S$ and $Y_S$ are independent,
\item $T_S$ is exponentially distributed with parameter $\sum_{j=1}^J k_{j}$,
\item $Y_S$ takes the value $i \in \{1, \ldots, J\}$ with probability $\frac{k_{i}}{\sum_{j'=1}^Jk_{j'}}$,
\end{enumerate}
which are exactly the properties needed to model the exit event using a kMC model, see Section~\ref{sec:kMC}.

\begin{figure}[h]
\centering
\psfrag{S}{$S$}
\psfrag{S_1}{$S_1$}
\psfrag{S_2}{$S_2$}
\psfrag{S_3}{$S_3$}
\psfrag{S_4}{$S_4$}
\psfrag{z_1}{$z_1$}
\psfrag{z_2}{$z_2$}
\psfrag{z_3}{$z_3$}
\psfrag{z_4}{$z_4$}
\psfrag{z_5}{$z_5$}
\psfrag{x_1}{$x_1$}
  \includegraphics[height=5cm]{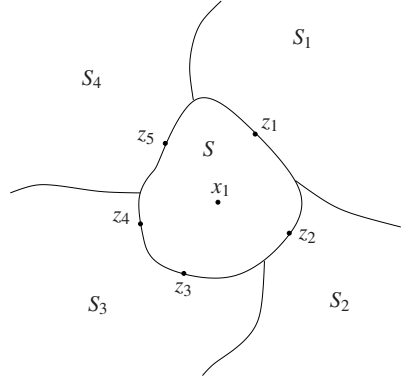}
  \caption{The points $(z_i)_{i=1, \ldots, 5}$ are the five local minima of $V$ on $\partial S$. The point $x_1$ is the global minimum of $V$ on $S$. In a global kMC model, exits through the neighborhood of $z_1$ is associated to a transition to $S_1$, exits through neighborhoods of $z_3$ and $z_4$ are associated to a transition to $S_3$, etc.}
  \label{fig:S_z}
\end{figure}

Notice that we have here considered a partition of the boundary $\partial S$ dictated by the {\em a priori} knowledge of neighboring states, having in mind a global kMC model as defined in Section~\ref{sec:kMC}. If one is only interested in simulating the exit from $S$ (without any {\em a priori} knowledge of $\R^d \setminus S$), a natural partition of the boundary is then 
$$\partial S=\cup_{i=1}^I B_{z_i}$$
where $(z_i)_{i=1, \ldots, I}$ are the local minima of $V$ on $\partial S$, see Figure~\ref{fig:S_z}, and for all $i \in \{1, \ldots ,I\}$, $B_{z_i}$ is a neighborhood of $z_i$. Indeed, in the small temperature regime, exits through $\partial S$ occur around the local minima of $V$ on $\partial S$ (this will be discussed in Section~\ref{sec:EK_math}, see Remark~\ref{rem:prefactors}). More precisely, for the mathematical analysis in Section~\ref{sec:EK_math}, we will define $B_{z_i}$ as the basin of
  attraction of $z_i$ for the dynamics $\dot{x} = -\nabla_T V(x)$ in the boundary $\partial S$ (where
  $\nabla_T V$ denotes the tangential gradient of $V$ along the
  boundary $\partial S$):
\begin{equation}\label{eq:Bzi}
B_{z_i}= \left\{ x_0 \in \partial S, \, \lim_{t \to \infty} x(t) = z_i, \text{ where } x(0)=x_0 \text{ and } \dot{x} = -\nabla_T V(x)\right\}.
\end{equation}
 We will then show that exits through $B_{z_i}$ actually only occur through a neighborhood of $z_i$ in the small temperature regime $\beta \to \infty$. Using this partition, one can thus simulate the exit event as explained above through a couple of random variable $(T_S,Y_S)$, where $Y_S=i$ if the exit occurs through the neighborhood $B_{z_i}$ of $z_i$ (where $i \in \{1, \ldots ,I\}$). The associated rates are those defined above, replacing the partition $(\partial S_i)_{i=1,\ldots,J}$ by $(B_{z_i})_{i=1,\ldots,I}$:
\begin{equation}\label{eq:ki_Bz}
k_i=\lambda_1 \P^{\nu_S}(X_{T_S} \in B_{z_i}) = - \lambda_1 
\frac{\displaystyle \int_{B_{z_i}} \partial_n u_1 \, d \sigma}{\displaystyle \beta \lambda_1
  \int_S u_1(x)  \, d x}.
\end{equation}

As explained above, this is useful to simulate the exit from the state $S$ much more efficiently than by integrating in time the Langevin~\eqref{eq:Langevin} or overdamped Langevin~\eqref{eq:overdamped_Lang}
dynamics. It can be used as such in a kinetic Monte Carlo model. It is also the basic ingredient of the accelerated dynamics that will be discussed in Section~\ref{sec:AD}. These algorithms aim at accelerating the sampling of the metastable trajectories of the Langevin or overdamped Langevin
dynamics, by efficiently generating the exit event from a metastable state $S$ when those dynamics get trapped in $S$. 

Notice that the results presented here are very general: they hold for both energetic and entropic traps, and for any Markov dynamics as soon as a QSD exists (non reversible Markov dynamics, discrete time Markov dynamics, etc).

One practical difficulty to use the results above as such in a kMC dynamics and in some of the accelerated dynamics algorithms, is that one needs in addition a simple way to evaluate the rates~$k_i$ (with a more explicit formula than the integral formulation in~\eqref{eq:ki_Bz}). We will come back to this in Section~\ref{sec:EK_HTST}, where the Eyring-Kramers formulas will be introduced to approximate these rates, in the small temperature regime.

\subsection{Estimating the convergence time to the QSD in practice}\label{sec:tau_corr}

We have seen in the previous section that if the process $(X_t)_{t \ge 0}$ starts from the QSD in $S$, one can model the exit event from $S$ exactly using a kMC model. In addition, in view of Proposition~\ref{prop:QSD_1}, the QSD is reached for the process conditioned to stay in $S$ in the longtime limit. A natural question is then: how long the process $(X_t)_{t \ge 0}$ should remain in $S$ so that one can assume it is sufficiently close to the QSD? In the original paper~\citep{voter-98}, this time is called the correlation time, denoted by $\tau_{corr}$, and is typically fixed as a constant for all the states, this constant being estimated using some physical intuition on the system at hand or an harmonic approximation. When the dynamics enters a state $S$, one thus waits for some time $\tau_{corr}$ before assuming that the QSD has been reached: this is called the decorrelation step in~\citet{voter-98}. Let us now discuss two ways to estimate the time needed to reach the QSD: one theoretical, and one numerical.

From a theoretical viewpoint, the following result proven in~\citet{le-bris-lelievre-luskin-perez-12} gives a first estimate of $\tau_{corr}$.
\begin{proposition}\label{prop:CV_QSD}
Let $(X_t)_{t \ge 0}$ satisfies~\eqref{eq:overdamped_Lang} with $X_0
\in S$. Let us consider $-\lambda_2 < -\lambda_1 < 0$  the first  two
eigenvalues of the operator $L^*$ on $S$ with homogeneous Dirichlet boundary
conditions on $\partial S$ (see Proposition~\ref{prop:QSD_2} for the
definition of $L^*$). Then, there
exists a constant $C>0$ which depends on the law of $X_0$, such that,
for all $t \ge \frac{C}{\lambda_2  - \lambda_1}$,
\begin{equation}\label{eq:CV_QSD}
\sup_{f,\, \|f\|_{L^\infty} \le 1} \left| \E(f(T_S-t,X_{T_S}) | T_S
  \ge t) -
  \E^{\nu_S}(f(T_S,X_{T_S})) \right| \le C \exp (-(\lambda_2-\lambda_1) t ).
\end{equation}
\end{proposition}
In other words, the total variation norm between the law of
$(T_S-t,X_{T_S})$ conditioned to $T_S \ge t$ (for any initial
condition $X_0 \in S$), and the law of
$(T_S,X_{T_S})$ when $X_0$ is distributed according to the QSD $\nu_S$,
decreases exponentially fast with rate $\lambda_2-\lambda_1$. This
means that $\tau_{corr}$ should be chosen of the order
$1/(\lambda_2-\lambda_1)$. There are however two difficulties with this result. First, this is not a very practical
result since computing the eigenvalues $\lambda_1$ and $\lambda_2$ is in general impossible. Second, the constant $C$ in the right-hand side of~\eqref{eq:CV_QSD} and in the inequality $t \ge \frac{C}{\lambda_2  - \lambda_1}$  depends in a complicated way on the initial condition $X_0 \in S$. And one indeed expects the convergence time to the QSD to strongly depend on the initial condition.

A more practical approach to estimate the convergence time $\tau_{corr}$ to the QSD has been proposed in~\citet{binder-simpson-lelievre-15}. The method uses two ingredients:
\begin{itemize}
\item The Fleming-Viot particle process~\citep{ferrari-maric-07}, which
  consists in
  $N$ replicas $(X^1_t, \ldots,X^N_t)_{t \ge 0}$ which are evolving and
  interacting in such a way that
  the empirical distribution $\frac{1}{N} \sum_{n=1}^N \delta_{X^n_t}$
  is close (in the large $N$ limit) to the law of the process $X_t$
  conditioned on $t<T_S$.
\item The Gelman-Rubin convergence diagnostic~\citep{gelman-rubin-92} to estimate the
  correlation time as the convergence time to a stationary state for
  the Fleming-Viot particle process.
\end{itemize}
Roughly speaking, the Fleming-Viot particle process consists in the following: each replica $(X^i_t)_{t \ge 0}$ evolves according to the
original dynamics~\eqref{eq:overdamped_Lang} driven by independent Brownian motions (starting from the initial conditions $X_0$), and, each time one of the replicas leaves the domain $S$, another one
taken at random is duplicated. The Gelman-Rubin convergence diagnostic
consists in comparing the average of a given observable over replicas
at a given time, with the average of this observable over time and
replicas: when the two averages are close (up to a tolerance, and for
a well chosen list of observables), the process is considered at stationarity. The interest of this approach is that it gives a practical way to approximate $\tau_{corr}$ which explicitly takes into account the initial condition of $(X_t)_{t \ge 0}$ in $S$. Of course, one should be cautious when using this technique in practice, since bad choices of the observables may lead to false convergence. We refer to~\citet{binder-simpson-lelievre-15} for examples of applications.
 
In the following, even though the correlation time $\tau_{corr}$ depends in principle on the state under consideration $S$ and on the initial condition of the process in $S$, we do not indicate explicitly this dependency and stick to the simple notation $\tau_{corr}$.

\section{Eyring-Kramers law and the Harmonic transition state theory}\label{sec:EK_HTST}

 In the previous section, we explained that if the process $(X_t)_{t \ge 0}$ remains sufficiently long in the state $S$, then the exit event can be exactly modeled using a kMC model, with the definitions~\eqref{eq:ki_Bz} for the rates associated with exits through local minima of $V$ on $\partial S$. As will become clear below, some accelerated algorithms will require in addition formulas for the rates, which explicitly depend on the potential $V$ and the temperature $\beta$. This is used in particular in  the Temperature Accelerated Dynamics (see Section~\ref{sec:TAD}) to infer the exit event at low temperature from exit events observed at a higher temperature. Such formulas are given by the so-called Eyring-Kramers laws, which are introduced in Section~\ref{sec:HTST}. A mathematical result showing that the exit rates from a state $S$ can indeed be approximated by the Eyring-Kramers laws will then be presented in Section~\ref{sec:EK_math}, for the overdamped Langevin dynamics~\eqref{eq:overdamped_Lang}.

\subsection{Eyring-Kramers laws}\label{sec:HTST}

The Eyring-Kramers laws give formulas which are used in many contexts as approximations of the rates modeling the exit event from the state $S$ in the small temperature regime ($\beta \to \infty$). Let us consider again the situation of Figure~\ref{fig:S_z}, and let us introduce the global minimum $x_1$ of $V$ in $S$, and the local minima $(z_i)_{i=1, \ldots, I}$ of $V$ on $\partial S$. We assume in the following that they are ordered such that
$$V(z_1) \le V(z_2) \le \ldots \le V(z_I).$$
Notice that we also assume in the following that all the critical points of $V$ and $V|_{\partial S}$ are non degenerate, which implies in particular that there are a finite number of local minima of $V|_{\partial S}$.
The Eyring-Kramers formula gives estimate for the exit rates $(k_{j})_{j=1, \ldots, I}$ through neighborhoods of the local minima $(z_{j})_{j=1, \ldots, I}$, namely: 
\begin{equation}\label{eq:EK}
\forall j \in \{1, \ldots, I\}, \, k_{j}=\nu_{j} \exp(-\beta[V(z_j)-V(x_1)])
\end{equation}
where $\nu_{j}>0$ is a prefactor which depends on the dynamic under
consideration and on $V$ around $x_1$ and $z_j$. Let us give a few
examples. If $S$ is taken as the basin of attraction of $x_1$ for the
dynamics $\dot{x}=-\nabla V(x)$ so that the points $z_j$ are order
one saddle points, the prefactor writes for the Langevin
dynamics~\eqref{eq:Langevin} (assuming $M={\rm Id}$ for simplicity):
\begin{equation}\label{eq:nu_L}
\nu^{L}_{j}=\frac{1}{4\pi} \left(\sqrt{\gamma^2 +4|\lambda^-(z_j)|}
  - \gamma \right)
\frac{\displaystyle \sqrt{\det (\nabla^2 V)(x_1)}}{\displaystyle \sqrt{|\det (\nabla^2 V)(z_j)|}}
\end{equation}
where $\nabla^2 V$ is the Hessian of $V$ (which, we recall, is assumed to be non degenerate) and $\lambda^-(z_j)$ denotes
the negative eigenvalue of $\nabla^2 V(z_j)$. This
formula was derived in~\citet{kramers-40} in a
one-dimensional situation.
The equivalent formula for the overdamped Langevin
dynamics~\eqref{eq:overdamped_Lang} is:
\begin{equation}\label{eq:nu_OL}
\nu^{OL}_{j}=\frac{1}{2\pi}|\lambda^-(z_j)| \frac{\displaystyle \sqrt{\det (\nabla^2 V)(x_1)}}{\displaystyle\sqrt{|\det (\nabla^2 V)(z_j)|}}.
\end{equation}
Notice that $\lim_{\gamma \to \infty} \gamma \nu^L_{j}=
\nu^{OL}_{j}$, as expected from the rescaling in time used to go
from Langevin to overdamped Langevin (see Section~\ref{sec:intro}). The
formula~\eqref{eq:nu_OL} has again been obtained in~\citet{kramers-40}, but also by many authors previously,
see the exhaustive review of the literature reported
in~\citet{hanggi-talkner-barkovec-90}.

Using the  Eyring-Kramers rates~\eqref{eq:EK} in the kinetic Monte Carlo model introduced in Section~\ref{sec:kMC}, one can then model the exit event from $S$ as follows. The exit event from $S$ is modeled by the couple of random variables $(T,Y)$ where $T$ is the first exit time from $S$ and $Y=j$ if the process exits  from $S$ in a neighborhood of $z_j$. This couple has the following law:
\begin{enumerate}
\item $T$ and $Y$ are independent
\item $T$ is exponentially distributed with parameter 
\begin{equation}\label{eq:er_EK}
\sum_{j=1}^n k_{j} \sim \left(\sum_{j'=1}^{I_0} \nu_{j'} \right) \exp (-\beta[V(z_1)-V(x_1)])
\end{equation}
\item $Y$ takes the value $j \in \{1, \ldots n\}$ with probability 
\begin{equation}\label{eq:ep_EK}
\frac{k_{j}}{\sum_{j'=1}^n k_{j'}} \sim \frac{\nu_j}{\sum_{j'=1}^{I_0} \nu_{j'}} \exp(-\beta[V(z_j)-V(z_1)]).
\end{equation}
\end{enumerate}
where the equivalents $\sim$ are valid in the small temperature regime $\beta \to \infty$, and where $I_0 \in \{1, \ldots, I\}$ denotes the number of global minima of $V$ on $\partial S$:
$$V(z_1) = \ldots = V(z_{I_0}) <V(z_{I_0+1}) \le \ldots \le V(z_I).$$
The modelling of the exit event using a kMC model parameterized by the Eyring-Kramers laws is sometimes called the Harmonic Transition State Theory in the literature. We refer for example to~\citet{voter-05} for an introduction to this theory and relevant references.

As already explained in Section~\ref{sec:kMC}, using such a model rather than the Langevin or overdamped Langevin dynamics to model the exit event from $S$ is useful either to simulate the evolution of the dynamics over very long time using the kMC model, or to accelerate the sampling of metastable trajectories of the Langevin or overdamped Langevin dynamics, as will be explained in Section~\ref{sec:AD}. 

This raises the following theoretical question: are the Eyring-Kramers rates~\eqref{eq:EK} a valid approximation to model the exit event from~$S$ using a kMC model, for the Langevin or overdamped Langevin dynamics? To be more precise, let us consider the Langevin or overdamped Langevin dynamics~\eqref{eq:overdamped_Lang}. We have already seen in Section~\ref{sec:QSD_kMC} that if the stochastic process $(X_t)_{t \ge 0}$ remains for a sufficiently long time in $S$, then it is indeed valid to use a kMC model to simulate the exit event from~$S$, with associated rates defined using the eigenvalue-eigenfunction pair $(\lambda_1,u_1)$, see Equation~\eqref{eq:ki_Bz}. The question is then: can these rates defined by~\eqref{eq:ki_Bz} be approximated using the Eyring-Kramers laws~\eqref{eq:EK}?

\subsection{Approximating the exit rates by the Eyring-Kramers laws: a mathematical framework}\label{sec:EK_math}

Let us consider the dynamics~\eqref{eq:overdamped_Lang} with an initial condition
distributed according to the QSD $\nu_S$ in a domain $S$. We assume
the following:
\begin{itemize}
\item The domain $S$ is an open smooth bounded domain in $\R^d$.
\item The function $V:\overline{S} \to \R$ is a Morse function\footnote{A Morse function is a function such that, at any critical point, the Hessian is non singular. Notice that this is in particular required to define the prefactors in the Eyring-Kramers laws.} with a single critical
point $x_1$. Moreover, $x_1 \in S$ and $V(x_1)=\min_{\overline{S}} V$. 
\item The normal derivative $\partial_n V$ is strictly positive on $\partial
  S$, and $V|_{\partial S}$ is a Morse function with local minima
  reached at $z_1, \ldots , z_I$ with $$V(z_1) = \ldots =V(z_{I_0}) < V(z_{I_0+1}) \le \ldots \le
  V(z_I).$$
\item The height of the barrier is large compared to the saddle points
  heights discrepancies: $V(z_1)-V(x_1) > V(z_I)-V(z_1)$.
\item For all $i \in \{1, \ldots I\}$, consider $B_{z_i}
  \subset \partial S$ the basin of
  attraction for the dynamics in the boundary $\partial S$: $\dot{x} = -\nabla_T V(x)$ (where
  $\nabla_T V$ denotes the tangential gradient of $V$ along the
  boundary $\partial S$, see Equation~\eqref{eq:Bzi}). We assume that
\begin{equation}\label{eq:agmon}
\inf_{z \in B_{z_i}^c} d_a (z_i,z) > \max( V(z_I)-V(z_i) , V(z_i)-V(z_1) )
\end{equation}
where $B_{z_i}^c=\partial S \setminus B_{z_i}$.
\end{itemize}
Here, $d_a$ is the Agmon distance:
$$d_a(x,y) = \inf_{\gamma \in \Gamma_{x,y}} \int_0^1 g(\gamma(t)) |\gamma'(t)| \, dt$$
where $g=\left\{ \begin{aligned} &|\nabla V| \text { in $S$} \\ &|\nabla_T
    V|  \text{ in $\partial S$} \end{aligned} \right.$, and the
infimum is taken over the set $\Gamma_{x,y}$ of all piecewise $C^1$ paths
$\gamma:[0,1] \to \overline{S}$ such that $\gamma(0)=x$ and
$\gamma(1)=y$.
The Agmon distance is useful in order to measure the decay of
eigenfunctions away from critical points. These are the so-called
semi-classical Agmon estimates,
see~\citet{simon-84,helffer-sjostrand-84}. Under the assumptions stated above, the following result is proven in \citet{di-gesu-le-peutrec-lelievre-nectoux-17}.

\begin{theorem}\label{th:EK}
Under the assumptions stated above,  in the limit $\beta \to \infty$,
the exit rate is 
\begin{equation}\label{eq:lambda1}
\lambda_1= \sqrt{\frac{\beta  \det  (\nabla^2V)   (x_1)}{2\pi}}
\sum_{k=1}^{I_0}\frac{  \partial_nV(z_k)
}{    \displaystyle\sqrt{ \det (\nabla^2 V_{|\partial S})   (z_k) }
}{\rm e}^{-\beta(V(z_1)-V(x_1))}( 1+ O(\beta^{-1}) ).
\end{equation}
Moreover, for any open set $\Sigma_i $ containing $z_i$ such that $\overline\Sigma_i \subset B_{z_i}$, 
\begin{align}
\P^{\nu_S} (X_{T_S} \in \Sigma_i)&=   
\frac{\frac{\partial_n V(z_i)} {\sqrt{ \det (\nabla^2
    V|_{\partial S })   (z_i) }}}{ \sum_{k=1}^{I_0} \frac {\partial_n V(z_k) } { \sqrt{ \det (\nabla^2 V|_{\partial S })
    (z_k) } } }  {\rm e}^{-\beta(V(z_i)-V(z_1))} (   1+   
 O(\beta^{-1})  ). \label{eq:EK-p}
\end{align}
\end{theorem}
 Formula~\eqref{eq:lambda1} (resp.~\eqref{eq:EK-p}) should be compared with~\eqref{eq:er_EK} (resp.~\eqref{eq:ep_EK}) introduced above. 
We refer to~\citet{di-gesu-le-peutrec-lelievre-nectoux-17} for a proof, and for other related results  (see also~\cite{helffer-nier-06} for a proof of~\eqref{eq:lambda1}).  The proof uses techniques developed in particular in the
previous works: \citet{helffer-klein-nier-04,helffer-nier-06,le-peutrec-10,lelievre-nier-15}. The analysis requires to combine various
tools from semiclassical analysis to address new questions: sharp
estimates on quasimodes far from the critical points for Witten
Laplacians on manifolds with boundary, a precise analysis of the normal
derivative on the boundary of the first eigenfunction of Witten
Laplacians, and fine properties of the Agmon distance on manifolds with
boundary.

Using the two results of Theorem~\ref{th:EK} and the formula~\eqref{eq:ki_Bz} for the definition of the exit rates, one obtains that the exit rate associated with an exit in the neighborhood $\Sigma_i$ of $z_i$ is
\begin{align}
k_{i}&=\lambda_1  \P^{\nu_S} (X_{T_S} \in \Sigma_i) \nonumber \\
&=\widetilde{\nu}^{OL}_{i} {\rm e}^{-\beta(V(z_i)-V(x_1))} (   1+   
 O(\beta^{-1})  ) \label{eq:k0i}
\end{align}
where the prefactors $\widetilde{\nu}^{OL}_{i}$ is given by
\begin{equation}\label{eq:nu_OL_tilde}
\widetilde{\nu}^{OL}_{i}=\sqrt{\frac{\beta}{2\pi}}  \partial_nV(z_i) \frac{ \sqrt{\det (\nabla^2 V)   (x_1) } } {  \sqrt{ \det
     (\nabla^2 V_{|\partial S})   (z_i) }  }.
\end{equation}
This should be compared to the formulas~\eqref{eq:EK}-\eqref{eq:nu_OL} of the previous section (see Remark~\ref{rem:prefactor} below for a discussion on the precise values of the prefactors).
This gives a rigorous framework to use a kMC model parameterized using Eyring-Kramers laws to model the exit event from $S$, as introduced in the previous section. Let us emphasize that Theorem~\ref{th:EK} provides estimates on the probability $\P^{\nu_S} (X_{T_S} \in \Sigma_i)$, and not only on the logarithm of these probabilities (as obtained for example using large deviation results, see Equation~\eqref{eq:LD} below). Notice that Theorem~\ref{th:EK} also give error estimates (actually, it can be shown that the terms $O(\beta^{-1})$ in~\eqref{eq:lambda1}, \eqref{eq:EK-p} and~\eqref{eq:k0i} admit full expansions in positive powers of $\beta^{-1}$). 

Let us finish this section with a few important remarks.

\begin{remark}[From generalized saddle points to real saddle points]\label{rem:prefactor}
As stated in the assumptions, Theorem~\ref{th:EK} is proven under the assumption that $\partial_n V >
0$ on $\partial S$: the local minima $z_1, \ldots ,z_I$ of $V$ on
$\partial S$ are therefore
not saddle points of $V$ but so-called {\em generalized saddle points} (see~\citet{helffer-nier-06,le-peutrec-10}). In a future work, we intend to extend these results to the case where
the points $(z_i)_{1 \le i \le I}$ are saddle points\footnote{This is indeed natural if $S$ is taken as the basin of attraction of $x_1$ for the gradient dynamics $\dot{x}=-\nabla V(x)$, as in Section~\ref{sec:HTST}.} of~$V$, in which
case we expect to prove the same result~\eqref{eq:k0i} for the exit
rates, with the prefactor $\widetilde{\nu}^{OL}_{i}$ being
$\displaystyle\frac{1}{\pi}|\lambda^-(z_j)| \frac{\displaystyle \sqrt{\det
    (\nabla^2 V)(x_1)}}{\displaystyle\sqrt{|\det (\nabla^2
    V)(z_j)|}}$ (this formula can be obtained using formal expansions on
the exit time and the Laplace's method). Notice that the latter formula differs from~\eqref{eq:nu_OL} by a
multiplicative factor $1/2$. This is due to the fact that $\lambda_1$ is the exit rate
from $S$ and not the transition rate to one of the neighboring state (see the
remark on page 408 in~\citet{bovier-eckhoff-gayrard-klein-04} on this multiplicative factor $1/2$ and the results on
asymptotic exit times in~\citet{maier-stein-93} for example). One way to understand this multiplicative factor $1/2$ is that once on the saddle point, in the limit $\beta \to \infty$, the process
has a probability one half to go back to $S$, and a probability one half
to effectively leave $S$. This
multiplicative factor does not have any influence on the law of the
next visited state which only involves a ratio of the rates $k_{i}$,
see Equation~\eqref{eq:Y_n}.
\end{remark}

\begin{remark}[On the importance of prefactors]\label{rem:prefactors}
The importance of obtaining a result including the prefactors in the
rates is illustrated by the following result, which is also proven
in~\citet{di-gesu-le-peutrec-lelievre-nectoux-17}.  Consider a simple
situation with only two local
minima $z_1$ and $z_2$ on the boundary $\partial S$, with $V(z_1) < V(z_2)$. Compare the two exit probabilities:
\begin{itemize}
\item The probability to leave
through $\Sigma_2$ such that $\overline{\Sigma_2} \subset B_{z_2}$
and  $z_2 \in \Sigma_2$;
\item The probability to leave through $\Sigma$ such that
  $\overline{\Sigma} \subset B_{z_1}$ and
$\inf_{\Sigma} V=V(z_2)$.
\end{itemize}
By classical results from the large deviation theory (see for
example~\eqref{eq:LD} below) the probability to exit through $\Sigma$ and
$\Sigma_2$ both scale like a prefactor times ${\rm
  e}^{-\beta(V(z_2)-V(z_1))}$:  the difference can only be read from the
prefactors. Actually, it can be proven that, in the limit $\beta \to\infty$,
$$\frac{\P^{\nu_S}(X_{T_S} \in \Sigma)}{\P^{\nu_S}(X_{T_S} \in \Sigma_2)} = O (
\beta^{-1/2} ).
$$
The probability to leave through $\Sigma_2$ (namely through the generalized saddle point~$z_2$)
is thus much larger than through $\Sigma$, even though the two regions
are at the same potential height.
This result explains why the local minima
of $V$ on the boundary (namely in our setting the generalized saddle points) play
an important role when studying the exit event.
\end{remark}

\begin{remark}[On the geometric assumption~\eqref{eq:agmon}]
Among the assumptions required to prove Theorem~\ref{th:EK}, the geometric assumption~\eqref{eq:agmon} involving the Agmon distances is probably the most unexpected one. Such an assumption indeed never appeared before in any mathematical works on Eyring-Kramers laws or on the analysis of the exit event. In~\citet{di-gesu-le-peutrec-lelievre-nectoux-17}, we investigated numerically this assumption, and we observed in some simple geometric settings that if it is not satisfied, the results of Theorem~\ref{th:EK} indeed do not hold (more precisely, by computing the exit probabilities~\eqref{eq:EK-p}, we observed that the prefactors are not those predicted by the Eyring-Kramers formulas).
\end{remark}

\begin{remark}[From overdamped Langevin to Langevin]\label{rem:lang_2}
The results of Theorem~\ref{th:EK} have been obtained for the overdamped Langevin dynamics~\eqref{eq:overdamped_Lang}, and it is then natural to ask whether similar results could be obtained for the Langevin dynamics~\eqref{eq:Langevin}. Indeed such estimates on the rates are actually assumed to be correct for Langevin dynamics in many models and numerical methods using Harmonic Transition State Theory, and in particular in the Temperature Accelerated Dynamics algorithm (see Section~\ref{sec:TAD}). We already mentioned in Remark~\ref{rem:lang_1} that the existence of the QSD for Langevin requires some additional investigations compared to overdamped Langevin. This is all the more true for the analysis of the exit event starting from the QSD in the small temperature regime. There is hope to obtain similar properties by combining results on the operator associated with Langevin dynamics with absorbing boundary conditions~\citep{nier-18}, and works on the semi-classical analysis of Langevin dynamics~\citep{herau-hitrick-sjostrand-11}.
\end{remark}

\begin{remark}[On the mathematical results on the Eyring-Kramers law]\label{rem:math_EK}
Given the importance of the Eyring-Kramers laws in the physics literature, many mathematical approaches have been proposed in order to justify these formulas.  Some authors adopt a global approach: they look at the
spectrum of the infinitesimal generator of the overdamped Langevin dynamics in the small temperature regime $\beta \to \infty$, see for example the
work by~\citet{helffer-klein-nier-04} based on semi-classical analysis
results for Witten Laplacian and the
articles by~\citet{bovier-eckhoff-gayrard-klein-04,bovier-gayrard-klein-05,eckhoff-05}  where a potential theoretic approach is adopted. These global approaches give the
 cascade of relevant time scales to reach from a local minimum any other local minimum which is lower in energy. However, they do not give the law of the exit event from a metastable state. 

In the context of this article, we are more interested in a local approach: we  consider the exit event from a state $S \subset \R^d$, and we would like to relate continous state space Markov dynamics such as~\eqref{eq:overdamped_Lang} and~\eqref{eq:Langevin} and kinetic Monte Carlo model to describe this exit event. In the mathematical literature, the most
famous approach to study the exit event is the large deviation
theory~\citep{freidlin-wentzell-84}.  In the small
temperature regime, large deviation results
provide the exponential rates~\eqref{eq:EK}, but without the
prefactors and without precise error bounds. For the
dynamics~\eqref{eq:overdamped_Lang}, a typical result on the exit point distribution
is the following (see~\citet[Theorem 5.1]{freidlin-wentzell-84}): for
all $S'$ compactly embedded in $S$, for any $\gamma >0$, for any $\delta > 0$, there
exists $\delta_0 \in (0,\delta]$ and $\beta_0 > 0$ such that for
all $\beta > \beta_0$, for all $x \in S'$ such that $V(x) <
\min_{\partial S} V$, and for all $y
\in \partial S$,
\begin{equation}\label{eq:LD}
{\rm e}^{-\beta ( V(y)-V(z_1) + \gamma)}\le \P^{x} (X_{T_S} \in
{\mathcal V}_{\delta_0}(y))  \le
{\rm e}^{-\beta (V(y)-V(z_1)- \gamma)}
\end{equation}
where ${\mathcal V}_{\delta_0}(y)$ is a $\delta_0$-neighborhood of $y$ in
$\partial S$. The strength of large deviation theory is that it is
very general: it applies to any dynamics (reversible or non
reversible) and in a very general geometric setting, even though it
may be difficult in such general cases to make explicit the rate
functional, and thus to determine the exit probabilities. 

Theorem~\ref{th:EK} thus differs from all these previous results in the following ways: it gives precise asymptotic estimates of the
distribution of $X_{T_S}$ ; the approach is local, justifies the
Eyring-Kramers formula~\eqref{eq:EK} with the prefactors and provides sharp error
estimates. 
\end{remark}

\section{Accelerated dynamics algorithms}\label{sec:AD}
The aim of this section is to explain how the mathematical analysis of metastability and metastable states presented in Sections~\ref{sec:kMC_QSD} and~\ref{sec:EK_HTST} can be used to build efficient algorithms to sample metastable dynamics such as~\eqref{eq:Langevin} or~\eqref{eq:overdamped_Lang}. We present for the sake of simplicity the algorithms in
the setting of the overdamped Langevin dynamics~\eqref{eq:overdamped_Lang}. The
algorithms can be generalized to the Langevin
dynamics~\eqref{eq:Langevin}, and actually to many Markov dynamics, under some assumptions that will be discussed below. We will present four algorithms together with their mathematical foundations: Parallel Replica (Section~\ref{sec:ParRep}), Parallel Trajectory Splicing (Section~\ref{sec:ParSplice}), Hyperdynamics (Section~\ref{sec:Hyper}) and Temperature Accelerated Dynamics (Section~\ref{sec:TAD}). Finally, we will discuss in Section~\ref{sec:metastable} how to define in practice the metastable states.

\subsection{Parallel Replica}\label{sec:ParRep}

The idea of the Parallel Replica algorithm is to evolve a reference replica following~\eqref{eq:overdamped_Lang}, and if it remains trapped for a long time in a state, to simulate in parallel the exit event from this state.
The original parallel replica algorithm (see~\citet{voter-98}) consists in iterating three steps:
\begin{itemize}
\item {\em The decorrelation step}: As already explained in Section~\ref{sec:tau_corr}, in this step, a reference replica
  evolves according to the original dynamics~\eqref{eq:overdamped_Lang}, until it
  remains trapped for a time $\tau_{corr}$ in one of the states. During this step, no error is made, since the reference
  replica evolves following the original
  dynamics (and there is of course no computational gain compared to a
  naive direct numerical simulation).  Once the reference replica has been
  trapped in one of the states (that we denote generically by $S$ in
  the following two steps) for a time $\tau_{corr}$, the aim is to generate very
  efficiently the exit event. This is done in two steps.
\item {\em The dephasing step}: In this preparation step, $(N-1)$
  configurations are generated within $S$ (in addition to the one
  obtained form the reference replica) as follows. Starting from the
  position of the reference replica at the end of the decorrelation
  step, some trajectories are simulated in parallel for a time
  $\tau_{corr}$. For each trajectory, if it remains within $S$ over the
  time interval of length $\tau_{corr}$, then its end point is
  stored. Otherwise, the trajectory is discarded, and a new attempt to get
  a trajectory remaining in $S$ for a time $\tau_{corr}$ is made. This
  step is pure overhead. The objective is only to get $N$
  configurations in $S$ which will be used as initial conditions in
  the parallel step.
\item {\em The parallel step}: In the parallel step, $N$ replicas are
  evolved independently and in parallel, starting from the initial conditions generated
  in the dephasing step, following the original dynamics~\eqref{eq:overdamped_Lang}
  (with independent driving Brownian motions). This step ends as soon
  as one of the replica leaves $S$. Then, the simulation clock is
  updated by setting the residence time in the state $S$ to $N$ (the
  number of replicas) times the exit time of the first replica which
  left $S$. This replica now becomes the reference replica, and one
  goes back to the decorrelation step above.
\end{itemize}
The computational gain of this algorithm is in the parallel step,
which (as explained below) simulates the exit event in a wall clock
time $N$ times smaller  in average than the time the reference walker would have spent to leave $S$. This of course requires a parallel
architecture able to handle $N$ jobs in parallel\footnote{For a
  discussion on the parallel efficiency, communication and
  synchronization, we refer to the papers~\citet{voter-98,perez-uberuaga-voter-15,le-bris-lelievre-luskin-perez-12,binder-simpson-lelievre-15}.}. This algorithm can be seen as a way to
parallelize in time the simulation of the exit event, which is not trivial because of the sequential nature of time.

In view of the results presented in Section~\ref{sec:kMC_QSD}, the Parallel Replica is indeed a consistent algorithm (it generates in a statistically correct way the exit event from metastable states) for the following reasons. First, in view of the
first property~\eqref{eq:QSD_1} of the QSD, the decorrelation step is
simply a way to decide whether or not the reference replica remains
sufficiently long in one of the states so that it can be considered as
being distributed according to the QSD, see also Section~\ref{sec:tau_corr} for an analysis of the error introduced by choosing $\tau_{corr}$ too small (see in particular Proposition~\ref{prop:CV_QSD}), and a discussion on how to evaluate $\tau_{corr}$ on the fly. Second, by the same arguments, the dephasing step is simply a rejection algorithm to generate
many configurations in $S$ independently and identically distributed
with law the QSD $\nu_S$ in $S$. Finally, the parallel step generates an exit event which is exactly
the one that would have been obtained considering only one
replica. Indeed, up to the error quantified in Proposition~\ref{prop:CV_QSD}, all
the replica are i.i.d. with initial condition the QSD
$\nu_S$. Therefore, according to the third property of
the QSD stated in Proposition~\ref{prop:QSD_3} (see item 1), their exit times $(T^n_S)_{n \in \{1, \ldots N\}}$ are i.i.d. with law an
exponential distribution ($T^n_S$ being the exit time
of the $n$-th replica) so that
\begin{equation}\label{eq:Nexp}
N \min_{n \in \{1, \ldots, N\}} (T^n_S) \stackrel{\mathcal L}{=}
T^1_S.
\end{equation}
This explains why the exit time of the first replica which leaves $S$
needs to be multiplied by the number of replicas $N$. This also shows
why the parallel step gives a computational gain in terms of wall
clock: the time required to simulate the exit event is divided by $N$
compared to a direct numerical simulation.
Moreover, since starting from the QSD, the exit time and the exit
point are independent (see item 2 in Proposition~\ref{prop:QSD_3}), one also has
$$X^{N_0}_{T^{N_0}_S} \stackrel{\mathcal L}{=}
    X^1_{T^1_S},$$
where $(X^n_t)_{t \ge 0}$ is the $n$-th replica and $N_0=\arg\min_{n \in \{1, \ldots, N\}}(T^n_S)$ is the index
of the first replica which exits $S$. The exit
point of the first replica which exits $S$ is statistically the same as
the exit point of the reference walker. Finally, by the independence property of
exit time and exit point, one can actually combine the two former
results in a single equality in law on couples of random variables, which shows that the parallel step is
statistically exact:
$$\left(N \min_{n \in \{1, \ldots, N\}} (T^n_S) ,
X^{N_0}_{T^{N_0}_S}\right)
\stackrel{\mathcal L}{=} (T^1_S, X^1_{T^1_S}).$$

We presented the results in the context of the overdamped Langevin
dynamics~\eqref{eq:overdamped_Lang} for simplicity.
 The mathematical analysis actually shows that the parallel replica is a very versatile algorithm. In
particular, it can be applied to both energetic and entropic barriers, and
it does not assume a small temperature regime (in contrast with the
analysis we will present below for Hyperdynamics and Temperature Accelerated Dynamics). In addition, it can also be used for non equilibrium system (driven by non conservative forces), and for any Markovian dynamics as soon as a QSD exist. We refer for example to~\citet{wang-plechac-aristoff-16} for recent applications to continuous-time Markov chains.

The only errors introduced in the algorithm are related to the rate of
convergence to the QSD of the process conditioned to stay in the
state. The algorithm is efficient if the convergence time to the
QSD is small compared to the exit time (in other words, if the states
are metastable). The Equation~\eqref{eq:CV_QSD}
gives a way to quantify the error introduced by the whole
algorithm. In the limit $\tau_{corr} \to \infty$, the algorithm generates
exactly the correct exit event. This implies that the resulting dynamics is exact in terms of the laws on trajectories (and not only of the time marginals -master equation- or the stationary state for example). It gives the correct transition times and trajectories to go from one state to another for example. The price to pay is that the details of the dynamics within the states are lost.

\begin{remark}[Parallel Replica for discrete-time Markov process]
As a remark, let us notice that in practice, discrete-time processes
are used (since the Langevin or overdamped Langevin dynamics are discretized in time). Then, the exit
times are not exponentially but geometrically distributed. It is
however possible to generalize the formula~\eqref{eq:Nexp} to this
setting by using the following fact:  if $(\sigma_n)_{n \in \{1, \ldots
  N\}}$ are i.i.d. with geometric law, then
$$N \left( \min(\sigma_1,
\ldots, \sigma_N)-1 \right) + \min \left( n \in \{1, \ldots ,N \} , \, \sigma_n =  \min(\sigma_1,
\ldots, \sigma_N) \right) \stackrel{{\mathcal L}}{=}\sigma_1.$$ We refer
to~\citet{aristoff-lelievre-simpson-14} for more details. 
\end{remark}

\begin{remark}[Generalized Parallel Replica dynamics]
In view of the analysis above, the crucial parameter of the Parallel Replica algorithm is the correlation time $\tau_{corr}$. As explained in Section~\ref{sec:tau_corr},~\citet{binder-simpson-lelievre-15} propose an approach to estimate~$\tau_{corr}$ on the fly by combining the Fleming-Viot particle process to simulate the law of the process $X_t$
  conditioned on $t<T_S$ and the Gelman-Rubin convergence diagnostic to estimate the
the convergence time to a stationary state for
  the Fleming-Viot particle process, which thus gives an estimate of the correlation time.
Using this idea, a generalized parallel replica algorithm is introduced
in~\citet{binder-simpson-lelievre-15}, as follows.
 Each time the reference replica enters a new state, a
Fleming-Viot particle process is launched using $(N-1)$ replicas simulated in
parallel (with initial condition the configuration of the reference replica). Then the decorrelation step consists in the following: if
the reference replica leaves $S$ before the Fleming-Viot particle
process reaches stationarity, then a new decorrelation step starts
(and the replicas generated by the Fleming-Viot particle process are discarded);
if otherwise the Fleming-Viot particle
process reaches stationarity before the reference replica leaves $S$,
then one proceeds to the parallel step. Notice indeed that the final positions
of the replicas simulated by the Fleming-Viot particle process can be
used as initial conditions for the processes in the parallel
step, since they are (approximately) distributed according to the QSD. This procedure thus avoids the choice of a correlation time $\tau_{corr}$ {\em a
  priori}: it is in some sense estimated on the fly, depending on the state under consideration, and on the initial condition within this state. For
more details, discussions on the correlations included by the
Fleming-Viot particle process between the replicas, and numerical experiments
(in particular in cases with purely entropic barriers), we refer to~\citet{binder-simpson-lelievre-15}.
\end{remark}

\subsection{Parallel Trajectory Splicing}\label{sec:ParSplice}

Parallel Trajectory Splicing (abbreviated as ParSplice) is a variant of the Parallel Replica algorithm which has been introduced in~\citet{perez-cubuk-waterland-kaxiras-voter-15}. The idea of the algorithm is to generate in parallel many trajectory segments, which spent at least a time $\tau_{corr}$ in one state before the beginning of the segment, and ends in a state where it again spends at least a time $\tau_{corr}$ before the end of the segment. A trajectory segment ending in a state $S$ is then appended to a trajectory segment starting in the same state $S$. The analysis of this algorithm is very similar to the one described above for Parallel Replica: if $\tau_{corr}$ is sufficiently large, the trajectory segments start under the QSD in a state, and end in a state being again distributed according to the QSD. Therefore, the exit event from the ending state can be simulated by considering any segment which starts under the QSD within the same state. This justifies the fact that these trajectory segments can be appended to form a long molecular dynamics trajectory.

The main difficulty when implementing this technique is to manage the concurrent generation of the segments from many molecular dynamics instances. A database of segments is populated as the simulation proceeds, and this database is used at the same time to generate the molecular dynamics trajectory. This should be done carefully in order to introduce no bias. For example, segments generated from a given state should not be used in increasing order of their generation time, but using a first-in-first-out queue (otherwise short segments would be favored). In addition, in order to schedule the production of new segments to be added to the database, one uses some prediction of the future trajectory, based on the current knowledge of the visited states. For more details, we refer to~\citet{perez-cubuk-waterland-kaxiras-voter-15}.

\subsection{Hyperdynamics}\label{sec:Hyper}

As in Parallel Replica, let us assume that a reference replica $(X_t)_{t \ge 0}$ following the overdamped Langevin dynamics~\eqref{eq:overdamped_Lang} remains trapped for a time $\tau_{corr}$ is a metastable state $S$. If $\tau_{corr}$ is sufficiently large, one can assume that the process is distributed according to the QSD.
The principle
of the hyperdynamics algorithm is then to raise the potential inside the
state in order to accelerate the exit from $S$. The algorithm
thus requires a biasing potential $\delta V: S \to \R$, which
satisfies appropriate assumptions detailed below. The algorithm then
proceeds as follows:
\begin{itemize}
\item Equilibrate the dynamics on the biased potential $V+\delta V$,
  namely run the dynamics~\eqref{eq:overdamped_Lang} on the process
  $(X^{\delta V}_t)_{t \ge 0}$ over the biased
  potential conditionally to staying in the well, up to the time the
  random variable $X^{\delta V}_t$ has distribution close to the QSD
  $\nu_S^{\delta V}$ associated
  with the biased potential. This first step is a preparation step, which is pure
  overhead.  The end point $X^{\delta V}_t$ will be used as the
  initial condition for the next step.
\item Run the dynamics~\eqref{eq:overdamped_Lang} over the biased
  potential $V + \delta V$ up to the exit time $T^{\delta V}_S$ from
  the state $S$. The simulation clock is updated by adding the
  effective exit time $B \times T^{\delta V}_S$ where $B$ is the so-called
  boost factor defined by
\begin{equation}\label{eq:B}
B=\frac{1}{T^{\delta V}_S} \int_0^{T^{\delta V}_S} \exp(\beta \,
\delta V(X^{\delta V}_t)) \, dt.
\end{equation}
The exit point is then used as the starting point for a new
decorrelation step.
\end{itemize}

Roughly speaking, the assumptions required on $\delta V$ in the
original paper~\citep{voter-97} are twofold:
\begin{itemize}
\item $\delta V$ is sufficiently small so that the exit event from the
  state $S$ can be modeled by a kinetic
  Monte Carlo models parameterized by the Eyring-Kramers laws.
\item $\delta V$ is zero on (a neighborhood) of the boundary $\partial
  S$.
\end{itemize}
The derivation of the method relies on explicit formulas for the laws of
the exit time and exit point, using the harmonic transition state theory, as explained in Section~\ref{sec:EK_HTST}.

The algorithm we present here is actually slightly different from the
way it is introduced in the original paper~\citep{voter-97}. Indeed, in
the original version, the local equilibration steps (decorrelation
step and equilibration step on the biased potential) are omitted: it is assumed that the states are sufficiently metastable (for both the original potential and the biased potential)
so that these local equilibrations are immediate. It would be
interesting to check if the modifications we propose here improve the
accuracy of the method.

Let us now discuss the mathematical foundations of this technique, and
in particular, a way to understand the formula~\eqref{eq:B} for the
boost factor. We actually need to compare two exit events. The first
one is the exit event for the original process $X_t$ following the
dynamics~\eqref{eq:overdamped_Lang}, starting from the QSD $\nu_S$
associated with the state~$S$ and the dynamics with potential $V$. The
second one is the exit event for the process $X^{\delta V}_t$ following the
dynamics~\eqref{eq:overdamped_Lang} on the biased potential $V +
\delta V$, starting from the QSD $\nu_S^{\delta V}$
associated with the state~$S$ and the dynamics with potential $V+\delta
V$. Referring to Section~\ref{sec:EK_HTST}, one way to justify the algorithm is to use the fact that, both on the original potential $V$ and the biased potential $V + \delta V$, a kMC model with transition rates defined by~\eqref{eq:EK} with prefactors such as~\eqref{eq:nu_L},~\eqref{eq:nu_OL} or~\eqref{eq:nu_OL_tilde} can be used to model the exit event, see in particular Theorem~\ref{th:EK}. Indeed, if this is the case, by using the biasing potential $V+\delta V$, one easily checks that the ratio of exit rates $k_j/k_i$ (which gives the relative probability to exit through a neighborhood of $z_j$ compared to the probability to exit through a neighborhood of $z_i$) do not depend on $\delta V$ since, by assumption, $\delta V$ is zero on $\partial S$: this gives the consistency of the algorithm in terms of the distribution of the first exit point. Concerning the exit time, we know that it is exponentially distributed with parameter $\sum_{i=1}^I k_i$, with explicit dependency on $V$ (or $V+\delta V$ for the biased potential) in the small temperature regime thanks again to the formulas~\eqref{eq:EK} (and \eqref{eq:nu_L},~\eqref{eq:nu_OL} or~\eqref{eq:nu_OL_tilde}). For example, in the framework of Theorem~\ref{th:EK}, the parameter of the exponential law is
$\lambda_1$, where the dependency of $\lambda_1$ on the potential $V$ is explicitly given by formula~\eqref{eq:lambda1}. In particular, one easily checks that (indicating the dependency of $\lambda_1$ on the underlying potential in parenthesis) $$\frac{\lambda_{1}(V + \delta V)}{\lambda_{1}(V)}=\sqrt{\frac{\det  (\nabla^2(V+\delta V) )  (x_1)}{\det  (\nabla^2(V))(x_1)}} {\rm e}^{\beta \delta V(x_1)}( 1+ O(\beta^{-1}) )$$
(we again used the fact that $\delta V$ is assumed to be zero on $\partial S$). The right-hand side is equal, in the regime $\beta \to \infty$, to $\frac{\int_S \exp(-\beta V)}{\int_S \exp(-\beta (V+\delta V))}( 1+ O(\beta^{-1}) )$, which is indeed well approximated by the boost factor $B$ (in the limit of a sufficiently large residence time $T^{\delta V}_S$, using the ergodicity of the dynamics). This shows that Theorem~\ref{th:EK} applied to~$V$ and to~$V+\delta V$ justifies the use of hyperdynamics.

Let us mention that in  Theorem~\ref{th:EK}, the relative error on the law of the exit event scales like $O(\beta^{-1})$. It is actually possible to show that the relative error for hyperdynamics scales like $O(\exp(-c \beta))$ for some positive $c$, and to prove this result under less stringent geometric conditions than  Theorem~\ref{th:EK}, see~\citet{lelievre-nier-15}. Again, generalizing these results to the Langevin dynamics~\eqref{eq:Langevin} is an open problem that we are currently investigating (see Remark~\ref{rem:lang_2} above).

Notice that, contrary to the Parallel Replica method, the
hyperdynamics is limited
to energetic barriers and a small temperature regime, at least for our mathematical analysis. On the other hand, for very high energetic
barriers, hyperdynamics is in principle much more efficient: the Parallel Replica
method only divides the exit time by $N$ (the number of replicas),
while for deep wells, the boost factor $B$ is very large.

\begin{remark}[On the biasing potential]
A practical aspect we did not discuss so far is the
effective construction of the biasing potential $\delta V$. In the
original article~\citep{voter-97}, A.F. Voter proposes a technique based on the Hessian
$\nabla^2 V$. Alternatively, a well-known method in the context of
materials science is the bond-boost method introduced by~\citet{miron-fichthorn-03}. More recently, some authors proposed to build the biasing potentiel on the fly, by using adaptive biasing techniques, see for example~\citet{tiwary-parrinello-13,bal-neyts-15,dickson-17}. Let us also mention the recent variant called {\em local hyperdynamics}~\citep{kim-perez-voter-13} which does not enter the previous framework since the bias is actually a non conservative biasing force (which does not derive from a biasing potential). 
\end{remark}

\subsection{Temperature Accelerated Dynamics}\label{sec:TAD}

Let us finally introduce the Temperature Accelerated Dynamics (TAD), 
see~\citet{sorensen-voter-00}. Let us assume again that a reference replica $(X_t)_{t \ge 0}$ following the overdamped Langevin dynamics~\eqref{eq:overdamped_Lang} remains trapped for a time $\tau_{corr}$ is a metastable state $S$, so that one can assume that the process is distributed according to the QSD.
The principle of TAD is to increase the
temperature (namely increase $\beta^{-1}$
in~\eqref{eq:overdamped_Lang}) in order to accelerate the exit from
$S$. The algorithm consists in
\begin{itemize}
\item Simulating many exit events from $S$ at high temperature, starting
  from the QSD at high temperature,
\item Extrapolating the high temperature exit events to low
  temperature exit events using the Eyring-Kramers law~\eqref{eq:EK}.
\end{itemize}
As for the hyperdynamics algorithm, in
the original paper by~\citet{sorensen-voter-00}, no equilibration step is
used: it is assumed that the states are sufficiently metastable at
both high and low temperatures so that the convergence to the QSD is
immediate. In particular, in the simulations at high temperature, the replica is simply bounced back into the state $S$ when it leaves $S$ (and not resampled from the QSD in~$S$). Let us now describe more precisely how the extrapolation
procedure is made.

Let us consider the exit event from $S$, at a given temperature. Using the results of Section~\ref{sec:kMC_QSD} the exit event can be modeled using a kMC model with transition rates $(k_i)_{i=1,\ldots,I}$ through the local minima $(z_i)_{i =1,\ldots,I}$ of $V$ on $\partial S$. More precisely, let us consider the process evolving in $S$ which is resampled according to the QSD within~$S$ after each exit event. For $i \in \{1, \ldots, I\}$, let us denote by $\tau_i$ the first exit time through a neighborhood of $z_i$ for this process. One readily checks that $\tau_i$ is exponentially distributed with parameter $k_i$, and that the $\tau_i$'s are independent. Therefore,  using the notation of Section~\ref{sec:kMC}, the exit time $T_S$ and exit local minimum $Y_S$ ($Y_S=i$ if the exit point is in $B_{z_i}$) satisfies:
$$(T_S,Y_S) \stackrel{\mathcal L}{=} (\min(\tau_1, \ldots, \tau_I),\arg\min(\tau_1, \ldots, \tau_I)).$$

Now, using the results of Section~\ref{sec:EK_HTST} (see in particular Theorem~\ref{th:EK}), namely using the Eyring-Kramers law to parameterize the kMC model, one can compute, in the small temperature regime, the ratios $\frac{k^{hi}_i}{k^{lo}_i}$, where $k_i^{hi}$ (resp. $k_i^{lo}$) denotes the rate at high temperature $\beta^{hi}$ (resp. low temperature $\beta^{lo}$). For example, if one considers real saddle points on the boundary, with the prefactors~\eqref{eq:nu_L} (for Langevin) or~\eqref{eq:nu_OL} (for overdamped Langevin), one obtains:
\begin{equation}\label{eq:extrapolation}
\frac{k^{hi}_i}{k^{lo}_i} \simeq
\exp(-(\beta^{hi}-\beta^{lo}) (V(z_i) - V(x_1))).
\end{equation}
Likewise, in the framework of  Theorem~\ref{th:EK} (namely for generalized saddle points), one obtains
\begin{equation}\label{eq:extrapolation_generalized}
\frac{k^{hi}_i}{k^{lo}_i} \simeq \sqrt{\frac{\beta^{hi}}{\beta^{lo}}}
\exp(-(\beta^{hi}-\beta^{lo}) (V(z_i) - V(x_1))).
\end{equation}
Using these formulas, one can infer  the exit events at inverse temperature $\beta^{lo}$ from the exit events observed at inverse temperature $\beta^{hi}$, since:
\begin{equation}\label{eq:TAD_justif}
(\tau_1^{lo},\ldots,\tau_I^{lo})\stackrel{\mathcal
  L}{=}(\Theta^1\tau_1^{hi},\ldots,\Theta^I\tau_I^{hi})
\end{equation}
where
$$\Theta^i=\frac{k^{hi}_i}{k^{lo}_i}\simeq\exp(-(\beta^{hi}-\beta^{lo}) (V(x_i) - V(x_0)))
$$
is a multiplicative factor constructed from the ratio of the
rates. In the equality in law
in~\eqref{eq:TAD_justif} the random variables $\tau^{hi/lo}_i$ are, as
described above, exponential random variables with parameter
$k^{hi/lo}_i$. To have analytical formula for the correction factors
$\Theta_i$ and make the algorithm practical, the Eyring-Kramers law is
assumed to be exact and one uses in practice
$\Theta^i=\exp(-(\beta^{hi}-\beta^{lo}) (V(x_i) - V(x_0)))$ (or $\Theta^i=\sqrt{\beta^{hi}/\beta^{lo}}\exp(-(\beta^{hi}-\beta^{lo}) (V(x_i) - V(x_0)))$).

The TAD
algorithm thus consists in running the dynamics at high temperature,
observing the exit events through the saddle points on the boundary of
the state, and updating the exit time and exit region that would have
been observed at low temperature. More precisely, if, at a given time, exits through the
saddle points $\{s_1,\ldots,s_k\} \subset\{1, \ldots ,I\}$ have been
observed, one computes
$\min(\Theta^{s_1}\tau_{s_1}^{hi},\ldots,\Theta^{s_k}\tau_{s_k}^{hi})$ and
$\arg\min(\Theta^{s_1}\tau_{s_1}^{hi},\ldots,\Theta^{s_k}\tau_{s_k}^{hi})$
to get the corresponding exit time and exit region at low temperature. 

The interest of TAD compared to a brute force saddle point search is
that it is not required to observe exits through all the saddle points
in order to obtain a statistically correct exit event. Indeed, a stopping criterium is introduced to
stop the calculations at high temperature when the extrapolation procedure
will not modify anymore the low temperature exit event
(namely will not modify $\min
(\Theta^{s_1}\tau_{s_1}^{hi},\ldots,\Theta^{s_k}\tau_{s_k}^{hi})$,
$\{s_1,\ldots,s_k\} \subset\{1, \ldots ,I\}$ being the saddle points
discovered up to the time the stopping criterium is fulfilled). This stopping
criterium requires to provide some {\em a priori} knowledge, typically a lower
bound on the barriers $V(z_j)-V(x_1)$ ($j \in \{1, \ldots, I\}$), or a lower
bound on the prefactors $\nu_j$ in~\eqref{eq:EK} (see for example~\citet{aristoff-lelievre-14} for a discussion).
In some sense,
TAD can be seen as a clever saddle point search, with a rigorous way
to stop the searching procedure.

From the above discussion, the mathematical analysis of the TAD algorithm thus requires to prove that the exit event from $S$ can be modeled using a kMC model parameterized by the Eyring-Kramers formulas. This is exactly the content of Theorem~\ref{th:EK} for the overdamped Langevin dynamics~\eqref{eq:overdamped_Lang} with generalized saddle points on $\partial S$.  See also~\citet{aristoff-lelievre-14}, for the case of a one dimensional potential. The generalization of Theorem~\ref{th:EK} to real saddle points on $\partial S$ and to the Langevin dynamics~\eqref{eq:Langevin} is a work under progress.

Notice that, compared to the hyperdynamics, TAD is really  based on the Eyring-Kramers formulas, with relative error terms which scales like $1/\beta$, while for hyperdynamics, one can prove that the error are exponentially small in $\beta$. On the other hand, the interest of TAD 
compared to the hyperdynamics is that it does not require a
biasing potential, which may be complicated to build in some situation.

\subsection{How to choose the metastable states ?}\label{sec:metastable}

Let us finally make a few comments on how the metastable states can be defined in practice. As explained above, the overall efficiency of the accelerated dynamics methods depends on the choice of these metastable states: the algorithms indeed provide an acceleration after assuming that the QSD has been reached in a metastable state. The efficiency of the algorithms thus depends on the choice of the states, which should be metastable regions, so that the stochastic process generically reaches the local equilibrium (the QSD) before leaving the state. How to define the metastable states is a difficult question, very much related to the definition of good reaction coordinates (or good reduced degrees of freedom) for molecular dynamics.  Notice that choosing good metastable states also implies being able to estimate the correlation time $\tau_{corr}$ within each state either from some {\em a priori} knowledge, or some on the fly estimates, as already explained in Section~\ref{sec:tau_corr}.

In the original papers by~\citet{voter-97,voter-98,sorensen-voter-00,perez-cubuk-waterland-kaxiras-voter-15}, the states are defined as the basins of attraction of the gradient dynamics 
\begin{equation}\label{eq:grad_dyn}
\frac{dq}{dt} = -\nabla V(q).
\end{equation}
For almost every initial conditions, this dynamics converges to a local minimum of $V$: there are thus as many states as local minima of $V$, and the states define (up to a negligible set of points) a partition of the state space $\R^d$. This means that in practice, a steepest descent is performed at regular time intervals to identify the state in which the system is. One big advantage of this definition is that the states do not need to be defined {\em a priori}: they can be numbered as they are discovered by the dynamics, when new local minima of $V$ are identified. Formally, in such a case, one can introduce a map
$$\mathcal S: \R^d \to \N$$
which to a given position associates a label of a basin of attraction of~\eqref{eq:grad_dyn}. Accelerated dynamics then aim at efficiently simulating the so-called state-to-state dynamics $(\mathcal S(q_t))_{t \ge 0}$ or $(\mathcal S(X_t))_{t \ge 0}$.

For a system with more diffusive or entropic barriers (this is typically the case for biological applications), one could think of defining the states using relevant reaction coordinates (see for example~\citet{kum-dickson-stuart-uberuaga-voter-04} where the states are defined in terms of the molecular topology of the molecule of interest). This requires to identify the states before starting the simulation.

An important point to make is that since the numerical methods are local in nature, one does not need a partition of the state space to apply these techniques. Two situations can then be considered. One can first define metastable states $(S_i)_{i \ge 1}$ which are disjoint open subsets of $\R^d$ (they are called
milestones by~\citet{faradjian-elber-04}, target sets or core sets by~\citet{schuette-noe-lu-sarich-vanden-einjden-11}). Each time the process enters one of these states, one checks if the QSD is reached before leaving (this is the decorrelation step, see Section~\ref{sec:tau_corr}) and then Parallel Replica, Hyperdynamics or Temperature accelerated dynamics can be used to efficiently sample the exit event. Parallel Trajectory Splicing can also be applied in such a situation. These techniques thus do not require a partition. On the other hand, since one has to simulate the original dynamics outside $\cup_{i \ge 1} S_i$,  if the time spent outside $\cup_{i \ge 1} S_i$ is large, the algorithms are less efficient.

Another possibility is to introduce again an ensemble $(C_i)_{i \ge 1}$ of disjoint open subsets of $\R^d$, and to define the state $S_i$ as follows\footnote{The state is then unbounded, which raises some mathematical question about the existence and uniqueness of the QSD, which could be addressed under appropriate assumptions on the growth of the potential $V$ at infinity.}:
$$S_i = \R^d \setminus \cup_{j \neq i} C_j.$$
The dynamics then goes as follows: when the process enters one of the core set $C_i$ one considers the exit event from the associated state $S_i$, namely the next visited core set which is different from $C_i$. In terms of kMC dynamics, this corresponds to projecting the dynamics $(X_t)_{t\ge0}$ (or $(q_t)_{t\ge0}$) onto a discrete state-space dynamics obtained by considering the last visited core set as in~\citet{vanden-eijnden-venturoli-ciccotti-elber-08,schuette-noe-lu-sarich-vanden-einjden-11}. The interest of such an approach is that the time to reach the QSD starting from $\partial C_i$ in $S_i$ should be quite small if the $C_i$'s are well chosen. In such a setting, the process thus should decorrelate very quickly before exiting, and an approximation using a global kMC algorithm becomes relevant. However, one drawback is that the description of the underlying continuous state space dynamics may become poor, especially for very small core sets (since the information of the last visited core set may not be very informative about the actual state of the system).

In summary, one should keep in mind that there is possibly room for improvements in the way the metastable states are defined, compared to the original papers where a partition of the state space in basins of attraction of the gradient dynamics~\eqref{eq:grad_dyn} is considered, see also~\cite{perez-uberuaga-voter-15} for a recent discussion. This is particularly relevant since there are now techniques to approximate the convergence time to the QSD which can be applied for sets which are not basins of attractions of~\eqref{eq:grad_dyn} (see~\citet{binder-simpson-lelievre-15} and Section~\ref{sec:tau_corr}).

\section{Conclusion and perspectives}\label{sec:conc}

We presented a mathematical analysis of the accelerated dynamics proposed by A.F.~Voter and co-workers in order to efficiently generate exit events from metastable states. The analysis is based on the notion of quasi-stationary distribution (QSD), which gives a natural framework to prove that the exit event from a state $S$ for the Langevin or overdamped Langevin dynamics can be modeled by a kinetic Monte Carlo (kMC) model, if the QSD is reached before exiting from $S$. Moreover, under some assumptions, we reported on recent mathematical results which show that the Eyring-Kramers formula can be used to parameterize the kMC model.

From a theoretical viewpoint, we already mention above that the mathematical analysis proving that the Eyring-Kramers formula are good approximations of the exit rates is for the moment restricted to the overdamped Langevin dynamics~\eqref{eq:overdamped_Lang}, and to a domain $S$ such that $\partial_n V >0$ on $\partial S$, which prevents us from considering real saddle points $z_i$ on $\partial S$. The generalization of the results presented in Sections~\ref{sec:kMC_QSD} and~\ref{sec:EK_HTST} to the Langevin dynamics~\eqref{eq:Langevin} and real saddle points on $\partial S$ is a work in progress. It would also be interesting to investigate Langevin or overdamped Langevin dynamics with non gradient forces (non-equilibrium systems).

From a numerical viewpoint, the mathematical analysis shows the versatility and the limitations of the accelerated dynamics algorithms. The Parallel Replica and the Parallel Trajectory Splicing algorithms for example are very general, and can be applied in many contexts: general Markov dynamics, general definitions of states, etc. Likewise, the principle of the Temperature Accelerated Dynamics algorithm can be used as soon as  formulas approximating the exit rates are available (one could think for example of versions of TAD where the potential is modified, in the spirit of Hyperdynamics but without the assumption that the biasing potential is zero on~$\partial S$). Generally speaking, the investigation of the performance of these algorithms to new physical settings and general Markov dynamics is a very promising research direction.

\begin{acknowledgement}
This work is supported by the European Research Council under the
European Union's Seventh Framework Programme (FP/2007-2013) / ERC
Grant Agreement number 614492. Part of this work was completed during the long programs ``Large deviation theory in statistical physics: Recent advances and future challenges'' at the International Centre for Theoretical Sciences (Bangalore) and ``Complex High-Dimensional Energy Landscapes'' at the Institute for Pure and Applied Mathematics (UCLA). The author would like to thank ICTS and IPAM for their hospitality.
\end{acknowledgement}

\end{document}